\documentclass[11pt,a4paper]{article}

\usepackage{times}
\usepackage{natbib}
\usepackage{enumerate}
\usepackage{textcomp}
\usepackage[latin1]{inputenc}
\usepackage[english]{babel}
\usepackage{amsmath,amsfonts}
\usepackage{graphicx}
\newcommand{\plus}{+}
\newcommand{\var}{\mbox{Var}}

\newcommand{\beq}{\begin{equation}}
\newcommand{\enq}{\end{equation}}
\newcommand{\bea}{\begin{eqnarray}}
\newcommand{\ena}{\end{eqnarray}}
\newcommand{\beas}{\begin{eqnarray*}}
\newcommand{\enas}{\end{eqnarray*}}
\newcommand{\ignore}[1]{}

\usepackage{enumerate}

\newcounter{theorem}
\def\thetheorem{\arabic{theorem}}
\newenvironment{theorem}{\refstepcounter{theorem}
{\bf Theorem \thetheorem} \em }
\newcounter{lemma}
\def\thelemma{\arabic{lemma}}
\newenvironment{lemma}{\refstepcounter{lemma}
\bigskip \noindent {\bf Lemma \thelemma} \em }{\bigskip}

\newcounter{corollary}
\def\thecorollary{\arabic{corollary}}

\newcounter{proposition}
\def\theproposition{\arabic{proposition}}

\newcounter{example}
\def\thedefinition{\arabic{example}}

\newcounter{rem}
\def\thedefinition{\arabic{rem}}

\newcounter{definition}
\def\thedefinition{\arabic{definition}}

\begin{document}

\renewcommand{\baselinestretch}{.93}

\begin{center}
LARRY GOLDSTEIN
\,  -- \,\,
YOSEF RINOTT
\end{center}

\newcommand{\Reg}{^\mathsf{R}}

\begin{center}
{\LARGE  \bf A Permutation Test For Matching and its Asymptotic
Distribution
\footnote{\textit{Received September 2003 and revised November 2003}}
}
\end{center}

%62G10, 60F05

\vspace{23pt}
{\footnotesize
\noindent {\sc Summary -}

We consider a permutation method for testing whether observations
given in their natural pairing exhibit an unusual level of
similarity in situations where any two observations may be similar
at some unknown baseline level. Under a null hypotheses where
there is no distinguished pairing of the observations, a normal
approximation with explicit bounds and rates is presented for
determining approximate critical test levels.

}

\vspace*{12pt}


\section{Introduction} The work in this paper was motivated by
examples such as the following.

{\bf Example 1.}  Schiffman et. al (1978), with statistical
assistance by one of the authors$^\dagger$, studied the influence
of a doctor's prior probabilities of diseases on diagnosis.
Statistical thinking, which can be formalized in Bayesian terms,
suggests that given a set of symptoms, a doctor's diagnosis or
ranking of possible diagnoses should be influenced not only by the
symptoms, but also by the disease prevalence at the time of
diagnosis. Doctors' information on prevalence may come, for
example, from textbooks, articles, and personal experience. The
goal of the study was to verify the influence of {\it personal}
prior information or opinion on disease prevalence (henceforth
referred to as ``personal prior") on diagnosis and help determine
whether doctors need to be better educated to take prevalence into
account, or if providing them with information on prevalence at
the time of diagnosis is useful.

In this study each doctor in a sample produced first a ranking $X$
of the prevalence, or of the probability of various diseases from
a given list; such a ranking represents the doctor's personal
prior. A compatible medical scenario was then presented to all
doctors, and each one of them produced a ranked list $Y$ of
possible diagnoses from the same given list. Rank correlations
between $X$ and $Y$ for each doctor were then computed. To test
the hypotheses that a doctor's personal prior does not influence
his diagnostic rankings, a null hypotheses of zero correlation
between each doctor's $X$ and $Y$ is not appropriate. Even with no
such influence, one would expect that pairs of rankings would have
some nonzero baseline correlation due to the influence of other
factors like common medical knowledge. The null hypothesis of
interest that there is no influence of personal prior is complex
since the baseline correlation is unknown. The presence of an
unknown baseline correlation raises the question of how high the
within-doctor rank correlations need to be to reject the null
hypothesis and assert the claim that there is influence of
personal prior on diagnostic rankings.

Correlations are used here as a measure of similarity between
ranked lists. Henceforth we will talk about similarity in general,
and the approach applies to any measure of similarity or proximity
defined on the sample space.

The main focus of this paper is on examples of the following kind:

{\bf Example 2.} This example is somewhat artificial, but it is
simpler and can clarify the issue; it will also help in explaining
the example that follows it which is rather similar. Consider an
instructor who wants to know if students are copying from their
neighbors in a class where students take an exam while seated in
pairs. Given a measure of similarity between exams, we expect any
two exams to be similar even in the absence of copying. Common
knowledge that all students hopefully have would make their exams
similar to a certain, unknown degree. Therefore, we want to test
if the similarity between seated pairs is unusual (due to copying)
relative to some unknown baseline similarity. This example is
different from the first in that here a similarity score can be
computed for any pair of exams $X_i,X_j$, whereas in the first
example the correlations of interest are those between $X$ and
$Y$.

{\bf Example 3.} Situations similar to Example 2 arise naturally
in environmental and medical studies, where subjects in a given
study group are matched (paired) by certain common background of
interest, such as having lived in the same neighborhood during a
given period, having certain common medical conditions or having
certain variables in common (e.g., gender, age, weight, etc.). In
order to assess the influence of the background in question on a
given set of certain medical conditions (denoted by $X_i$ for
subject $i$), one should test whether matched pairs are more
similar than unmatched ones relative to the medical condition
being studied. The baseline similarity between unmatched pairs is
again unknown, but a certain degree of similarity must certainly
exist due to common factors that all subject in the particular
study might have. More specifically, suppose we have an even
number $n$ of subjects and those indexed by $2i-1$ and $2i$ form
the matched pairs for $i=1, \ldots, n/2$, and let $X_i$ measure
subject $i$'s medical condition. Our goal is to test whether all
$X_{2i-1}$ and $X_i$, which arise from the matched pairs, exhibit
more similarity then $X_i$ and $X_j$ from unmatched pairs.

A related testing problem arises in the design of studies
involving matching of subjects that are similar by some background
criteria in order to reduce variability in other variables of
interest. The matching process often requires a great effort. The
question of whether it achieves its purpose in producing a higher
level of similarity in the variables of interest than would be
achieved at random, can be tested as described in this paper.

Example 1 is a specific instance of a problem of the following
type. Consider pairs of observations $(X_1,Y_1),\ldots,
(X_n,Y_n)$, where $X_i$ and $Y_i$ take values in a space so that a
proximity function $c(X,Y)$ is defined. This function may
sometimes be obtained as a decreasing function of some metric.
However, for the rankings of Example 1, the rank correlation is a
relevant proximity function not derived from a metric. We want to
test whether the natural pairing of $X_i$ to $Y_i$ exhibits a
significantly higher level of proximity or similarity than an
unknown baseline level. The null hypotheses that the level of
proximity or similarity in the natural pairing is the same as the
baseline level can be formulated as the hypothesis that the
observations $[(X_1,Y_{\pi(1)}), \ldots,(X_n,Y_{\pi(n)})]$ are
identically distributed for all $\pi \in S_n$, the permutation
group of $n$ elements. Conditioning on the  observed
$\{e_{ij}=c(X_i,Y_j)\}$, a permutation test which compares the
value of \bea \label{U} U_\pi = \sum_{i=1}^n e_{i \pi(i)}\,\ena
for the special permutation $\pi=\mbox{id}$ (the identity),
against critical values of the distribution of $U_\pi$ when $\pi$
is uniform over $S_n$ (the distribution assigning equal
probabilities to every element in $S_n$), can be used to test the
null hypothesis.

The permutation distribution of $U_\pi$ for $\pi$ uniform over
$S_n$ was studied in numerous other statistical contexts. For a
seminal reference which contains both theory and applications see
Wald and Wolfowitz (1944). More recent articles which in turn
contain further references include the following. In connection
with linear rank statistics, Ho and Chen (1978) and Bolthausen
(1984) computed bounds on the rate of convergence to normality.
Bickel and van Zwet (1978) give more background and results on
linear rank statistics for two-sample problems, including an
Edgeworth expansion for a special case of (\ref{U}). Diaconis,
Graham and Holmes (1999) discuss similar statistics and also some
subsets of permutations related to tests for independence. Kolchin
and Chistyakov (1973) discuss the permutation distribution for the
subset of permutations with one cycle. Below we discuss a rather
different subset of permutations, in which the number of cycles is
maximal. For general theory on permutation tests see, e.g, Pesarin
(2001) and references therein. Related work on normal
approximations can be also be found in Stein (1986), whose ideas
and methods have strongly influenced us and other authors. We will
give a brief indication of some basic ideas of Stein's method.

In this paper we focus on situations as described in Examples 2
and 3 where all pairings can be compared. We consider the
following framework. Given an even number $n$ of paired
observations $(X_1,X_2),(X_3,X_4),\ldots,(X_{n-1},X_n)$, with
values in a space so that a proximity function $c(X_i,X_j)$ is
defined, we want to test whether the special pairing of $X_{2i-1}$
with $X_{2i}$ exhibits a significantly higher level of similarity
than an unknown baseline level. The null hypotheses that the
similarity level of the special pairing is the same as the
baseline level is here formulated as the hypothesis that the
variables $[(X_i,X_{\pi(i)}), i < \pi(i)]$ are identically
distributed for all $\pi \in \Pi_n$ where $$\Pi_n = \{ \pi \in
S_n: \pi^2=\mbox{id}, \pi(i) \not = i \quad \mbox{for all $i$}
\}.$$ The condition $\pi^2=\mbox{id}$ reflects the fact that if
$i$ is paired with $j$ then $j$ is paired with $i$, and the
condition $\pi(i) \not = i$ the fact that no $i$ can be paired
with itself. The special pairing which we suspect may show a high
similarity level corresponds to the permutation ${\tilde \pi} \in
\Pi_n$ specified by the conditions ${\tilde \pi}(2i-1)=2i$ and
$({\tilde \pi})^2=\mbox{id}$. Conditioning on the set of values
$\{e_{ij}=c(X_i,X_j)\}$ we consider the permutation test which
compares the value of $U_\pi$ at the special permutation
$\pi={\tilde \pi}$ against critical values of the distribution of
$U_\pi$ when $\pi$ is uniform over $\Pi_n$.

The proposed two tests discussed above appear similar, as in both
tests the observed similarity related to a special permutation is
compared to critical values computed against a null distribution
induced by the uniform distribution over a space of relevant
permutations. For the first test that space is $S_n$ and the
special permutation is the identity (which matches $X_i$ with
$Y_i$), and for the second test the space of permutations is
$\Pi_n$ and the special permutation is $\tilde{\pi}$ (which
matches $X_{2i-1}$ with $X_i$). We henceforth discuss only the
second case and study the permutation distribution relative to
$\Pi_n$. The methods used here apply to the permutation
distribution over the whole of $S_n$ mutatis mutandis.

For the null hypothesis to be true it is sufficient that the $X's$
are exchangeable, but the null hypothesis is complex and does not
specify the distribution of $U_\pi$ nor the baseline similarity.
In the absence of a null distribution, the above permutation test
seems very natural.

We shall provide a normal approximation to the permutation
distribution of $U_\pi$ of (\ref{U}) including bounds, rates, and
explicit constants in order to determine approximate critical
values for the permutation test.

Henceforth we suppress the dependence of $U_\pi$ in (\ref{U}) on
$\pi$. Furthermore, for values $g_{ij}$ with $g_{ii}=0$, we set $$
g_{i+}=\sum_{j=1}^n g_{ij}, \quad g_{+j}=\sum_{i=1}^n g_{ij},
\quad g_{++}=\sum_{i,j =1}^n g_{ij},\quad \mbox{and} \quad
\overline{g_{i+}} = \frac{1}{n-1}g_{i+}. $$ \ignore{For
consistency with later notation we write sums over both $i$ and
$j$ with $i \not = j$ as a sum over $|\{i,j\}|=2$.} Note that the
terms $e_{i\pi(i)}$ and $e_{\pi(i)i}$ always appear together in
the sum $U$, and we may therefore assume without loss of
generality that $e_{ij}=e_{ji}$. The diagonal terms $e_{ii}$ never
enter $U$ and we take them to be 0. Given such a collection of
numbers $e_{ij}$, define
\begin{equation}
\label{forgotten} d_{ij}= \left\{ \begin{array}{ccc}
e_{ij}-\frac{e_{i+}}{(n-2)}-\frac{e_{+j}}{(n-2)}+\frac{e_{++}}{(n-1)(n-2)}
& i
\not = j\\
0 & i = j. \end{array} \right.
\end{equation}

Bounds to the normal approximation for the permutation
distribution of $U$ are contained in the following theorem. For
convenience we assume without further comment that $n \ge 10$.

\begin{theorem}
\label{theorem-perm} Let $U$ be given by (\ref{U}), $\pi$ be
uniform over $\Pi_n$ and $$ W = \frac{U-EU}{\sqrt{\var(U)}},\quad
\alpha = \max |d_{ij}-d_{kl}|, \quad \mbox{and} \quad \delta =
\sup_{w \in R} | P(W \le w)- \Phi(w)|, $$ where $\Phi$ is the
standard normal distribution function. Then $$ EU=e_{++}/(n-1), $$
\begin{equation}
\label{VarUd} \mbox{Var}(U) = \frac{2}{(n-1)(n-3)} \left((n-2)
\sum_{i,k=1}^n e_{ik}^2 + \frac{1}{n-1}e_{\plus \plus}^2-
2\sum_{i=1}^n e_{i \plus}^2 \right),
\end{equation}
and there exist constants $c_1,c_2$ such that \bea \label{th1a}
\delta  &\le& c_1 n^{1/2} \sqrt{\{ \sum_{i,j=1}^n
d_{ij}\,^4/(\sum_{i,j=1}^n d\,^2_{ij})^2 \} } + \frac{c_2
\alpha^3\,n^{5/2}}{(\sum_{i,j=1}^n d\,^2_{ij})^{3/2}}. \nonumber
\ena
\end{theorem}

If, for example, the constants $d_{ij}$ are bounded then $\alpha$
is bounded and $\sum_{ij} d_{ij}^2=O(n^2)$, so in view of
(\ref{varU}), the bound above decays at the rate of
$\var(U)^{-1/2}=n^{-1/2}$. Below a somewhat crude calculation
gives the upper bounds of $c_1 \le 86, c_2 \le 243$.

\section{Proof of Theorem 1}
We compute the mean and variance of $U$ in Section
\ref{mean-variance}, and establish the upper bound on the normal
approximation in Section \ref{establish}.
\subsection{Mean and Variance of $U$}
\label{mean-variance} To compute the mean and variance of
$U=\sum_{i=1}^n e_{i \pi(i)}$, where $\pi$ is chosen uniformly
from $\Pi_n$, we have the following Lemma.

\begin{lemma}
\label{effigy} Let $g_{ij}$ satisfy $g_{ii}=0$ and set
$$ f_{ij}=
\left\{ \begin{array}{ccc} g_{ij}-\overline{g_{i+}} & i
\not = j\\
0 & i = j. \end{array} \right.
$$

Then with
$$
V = \sum_{i=1}^n g_{i \pi(i)}
$$
we have
$$
Eg_{i \pi(i)} = \overline{g_{i+}} \quad \mbox{and therefore} \quad
EV= \sum_{i=1}^n \overline{g_{i+}}
$$
and
$$
\mbox{Var}(V)= \frac{1}{(n-1)(n-3)} \left((2n-5)
\sum_{|\{i,j\}|=2}f_{ij}^2 + \sum_{|\{i,j\}|=2}f_{ij}f_{ji}
\right).
$$
\end{lemma}
\noindent {\bf Proof}: Since $\pi(i)$ can be any $j \not = i$ with
probability $1/(n-1)$, we have
$$ Eg_{i \pi(i)} = \frac{1}{n-1} \sum_{j:j \not = i} g_{ij} =
\frac{1}{n-1}\,g_{i+}= \overline{g_{i+}},
$$
and so \beas
\var(V) &=& \var \sum_{i=1}^n f_{i\pi(i)} \\
&=& \sum_{i=1}^n E f^2_{i\pi(i)} + \sum_{|\{i,j\}|=2}
E(f_{i\pi(i)}f_{j\pi(j)})\\
&=& \frac{1}{n-1} \sum_{|\{i,j\}|=2} f_{ij}^2 + \sum_{|\{i,j\}|=2}
E(f_{i\pi(i)}f_{j\pi(j)}). \enas

Now note that the probability is $1/(n-1)$ that $\pi(i)=j$, and
therefore that $\pi(j)=i$. If $\pi(i) \not = j$, then $i,j,
\pi(i), \pi(j)$ are all distinct, and given any $|\{i,j,k,l\}|=4$
the probability that $\pi(i)=k$ and $\pi(j)=l$ is
$1/[(n-1)(n-3)]$. We therefore have \bea \label{negcor}
E\sum_{|\{i,j\}|=2} f_{i\pi(i)}f_{j\pi(j)} =
\frac{1}{n-1}\sum_{|\{i,j\}|=2}f_{ij}^2 +
\frac{1}{(n-1)(n-3)}\sum_{|\{i,j,k,l\}|=4}f_{ik}f_{jl}. \ena

The first equality below follows by summing over $l \not \in
\{i,j,k\}$, and using $f_{jj}=0$ and $f_{j+}=0$, and the second in
a similar way by summing over $j \not \in \{i,k\}$; \beas
\sum_{|\{i,j,k,l\}|=4}f_{ik}f_{jl} &=&
\sum_{|\{i,j,k\}|=3}f_{ik}(-f_{ji}-f_{jk})\\ &=&
\sum_{|\{i,k\}|=2} f_{ik}(f_{ki}+f_{ik})\\ &=& \sum_{|\{i,k\}|=2}
\left( f_{ik}f_{ki} + f_{ik}^2 \right). \enas The formula for
$\mbox{Var}(V)$ now follows by collecting terms. $\Box$

Writing for the moment $U_d$ and $U_e$ for the values of $U$ based
on $d_{ij}$ and $e_{ij}$ respectively, we have
$$
U_d = U_e - \frac{e_{++}}{n-1}.
$$
In order to see the above relation between $U_d$ and $U_e$, sum
(\ref{forgotten}) over $i$ with $i \ne j$, and use symmetry to
yield $e_{+j}=e_{j+}$, and obtain
$$
d_{+j}=e_{+j} - [e_{++}-e_{j+}]/(n-2)-(n-1)e_{+j}/(n-2)+
e_{++}/(n-2) =0,
$$
so that
\begin{equation}
\label{mamamia} d_{i+}=d_{+j}=d_{++}=0.
\end{equation}
Since the distribution of $U_d$ is a simple translation of that
for $U_e$ we study $U_d=U$; henceforth we suppress the $d$.

Applying Lemma \ref{effigy} with $g_{ij}=d_{ij}$, since $d_{i+}=0$
we have $f_{ij}=d_{ij}$ and therefore $$ EU = 0; $$ using also
$d_{ij}=d_{ji}$, \begin{equation} \label{varU} \var (U) =
\frac{2(n-2)}{(n-1)(n-3)} \sum_{i,j=1}^n d_{ij}^2. \end{equation}

In terms of the symmetric but otherwise arbitrary values $e_{ij}$
which may not satisfy (\ref{mamamia}), the variance in
(\ref{VarUd}) is obtained by substituting (\ref{forgotten}) into
(\ref{varU}).

\subsection{Normal Approximation Upper Bound}
\label{establish} We apply the following theorem, which is a
special case of (1.10) of Theorem 1.2 of Rinott and Rotar (1997),
when $R=0$, using (1.12). The latter is based on Stein's method
(Stein 1986, pg 35), with an improvement on the rates under some
condition.

\begin{theorem}
\label{base} Let $(W,W^*)$ be exchangeable with $EW=0$ and
$EW^2=1$ such that for $0 < \lambda <1$ we have
\begin{equation}
\label{linearity} E(W^*|W) = (1-\lambda)W.
\end {equation}

If
\begin{equation}
\label{bou} |W^*-W|\,\leq\,A
\end{equation}
for a constant $A$, then \bea \label{th1} \lefteqn{\delta =
\sup_{w \in R} | P(W \le w)- \Phi(w) |} \nonumber\\ &\leq& \frac
{12}{\lambda} \sqrt{\var \{ E[(W^*-W)^2|W]\} } +
\sqrt{\frac{2}{\pi}}\frac{(48+\sqrt{32})A^3}{\lambda}. \ena
\end{theorem}

We briefly indicate the idea behind the proof of a theorem of this
type. This discussion can serve as some introduction to Stein's
method for the interested reader, but it is not necessary for the
rest of the paper.

First note that a random variable $W$ has the standard normal
distribution if and only if
\begin{equation}
\mbox{\rm E}f^{\prime }(W)=\mbox{\rm E}Wf(W)  \label{basic}
\end{equation}
holds for all continuous and piecewise continuously differentiable
functions $f$, for which the expectations in (\ref{basic}) exist.
This motivates the differential equation (\ref{DE}) in the lemma
below.

Set ${\bf \Phi }h= \mbox{\rm E}h(Z)$, where $Z$ is standard
normal, and $h$ is a function for which the expectation exists.
Also, for a real valued $h$, let $||h||$ denote the sup norm, that
is, $||h||=\sup_x|h(x)|.$ The lemma below is elementary, though
the bounds in (\ref{bounds}) require some calculations.

\begin{lemma}
\label{basic2} Let $h$ be a bounded piecewise continuously
differentiable real valued function. The function
\begin{equation}
f(w)=e^{w^{2}/2}\int_{-\infty }^{w}[h(x)-{\bf \Phi
}h]e^{-x^{2}/2}dx \label{DEsol}
\end{equation}
solves the (first order linear) differential equation
\begin{equation}
f^{\prime }(w)-wf(w)=h(w)-{\bf \Phi }h,  \label{DE}
\end{equation}
and
\begin{equation}
({\rm a})\,\quad ||f||\leq \sqrt{2\pi }||h||,\quad \,\,({\rm b}%
)\,\,\quad ||f^{\prime }||\leq 2||h||,\quad \,\,({\rm c})\quad
\,\,||f^{\prime \prime }||\leq 2||h^{\prime }||.
\label{bounds}
\end{equation}
\end{lemma}
Now, exchangeability of $(W,W^*)$ and (\ref{linearity}) directly
imply $$E \{Wf(W)\}=\frac{E \{(W^{* }-W)[f(W^{*
})-f(W)]\}}{2\lambda }.$$ Together with (\ref{DE}) this implies
\begin{equation}
E h(W)-{\bf \Phi }h=E f^{\prime }(W)-\frac{E \{(W^{* }-W)[f(W^{*
})-f(W)]\}}{2\lambda }.  \label{ex}
\end{equation}

%For $W^{* }$  close to $W$,
The first term in the Taylor expansion of $ f(W^{* })-f(W)$ is
$(W^{* }-W)f^{\prime}(W)$, and the r.h.s. of (\ref{ex}) is bounded
by $\frac{1}{2\lambda }E \{f^{\prime }(W)[2\lambda -(W^{*
}-W)^{2}]\}$ plus a remainder term which we now ignore. By
(\ref{linearity}) we have $E (W^{* }-W)^{2}=2\lambda ,$ and
(\ref{bounds} b) and the Cauchy Schwarz inequality readily yield
\[
E \{f^{\prime }(W)[2\lambda -(W^{* }-W)^{2}]\}\leq 2 ||h||
\sqrt{{\rm Var}\{E [(W^{* }-W)^{2}|W]\}}.
\]
Using (\ref{ex}), an approximation of an indicator function of a
half line by a smooth $h$ yields a term similar to the first term
on the r.h.s. of (\ref{th1}), and calculation of the remainder in
the above Taylor expansion yields the second. To obtain the
precise bound (\ref{th1}), a certain induction and further
calculations are needed. $\Box$

We shall apply Theorem \ref{base} to $W=U/\sigma$, but for
convenience we first describe the coupling and compute the
relevant quantities in terms of $U$. Given a permutation $\pi$
chosen uniformly from $\Pi_n$ construct the permutation $\pi^*$ in
$\Pi_n$ by choosing $I,J$ distinct and uniformly, and imposing
$\pi^*(I)=J$ (and therefore $\pi^*(J)=I$), and
$\pi^*(\pi(I))=\pi(J)$ (and therefore $\pi^*(\pi(J))=\pi(I)$) and
fixing the values of $\pi^*(k)=\pi(k)$ for $k \not \in
\{I,J,\pi(I),\pi(J)\}$. With $U = \sum_i d_{i \pi(i)}$, let $U^*=
\sum_i d_{i \pi^*(i)}$.

To verify (\ref{linearity}), first note that \bea \label{wsmw} U^*
- U &=& 2\left( d_{IJ} +
d_{\pi(I)\pi(J)}-(d_{I\pi(I)}+d_{J\pi(J)})\right), \ena where the
factor 2 accounts for the symmetry $d_{ij}=d_{ji}$.
%For the last two terms we have \bea \label{no-indicator}
%E(d_{I\pi(I)}|U)=E(d_{J\pi(J)}|U) = \frac{1}{n}U. \ena

Letting $C$ be the event that $J \not = \pi(I)$, we have
$(U^*-U)=(U^*-U){\bf 1}_C$ and therefore
$$
E((U^*-U)|U) = E((U^*-U){\bf 1}_C|U).
$$

For the first two terms in (\ref{wsmw}), recalling $d_{++}=0$, and
using that $(I,J)$ is independent of $\pi$ and equals any of the
$n(n-1)$ pairs $(i,j)$ for which $i \not = j$, \beas E(d_{IJ}{\bf
1}_C|\pi)\ &=&
\frac{1}{n(n-1)} \sum_{|\{i,j\}|=2}d_{ij}{\bf 1}(j \not = \pi(i))\\
&=&\frac{1}{n(n-1)} \sum_{i=1}^n \left( \sum_{j: j \not = i}
d_{ij}-d_{i \pi(i)} \right) \\
&=& \frac{-1}{n(n-1)} U,\enas and similarly for the term
$d_{\pi(I)\,\pi(J)}$, as $(\pi(I),\pi(J))$ has the same
distribution as $(I,J)$.

Now consider the third term on the right hand side in
(\ref{wsmw}): \beas && E(d_{I\pi(I)}{\bf 1}_C|\pi) =
\frac{1}{n(n-1)}\sum_{|\{i,j\}|=2}
d_{i\pi(i)} {\bf 1}(j \not = \pi(i))\\
&=& \frac{1}{n(n-1)}\sum_{i=1}^n d_{i\pi(i)} \sum_{j: j \not =
i}{\bf 1}(j \not = \pi(i)) = \frac{1}{n(n-1)}\sum_{i=1}^n
d_{i\pi(i)} \sum_{j: j \not \in \{ i, \pi(i)\} }1 \\
&=& \frac{n-2}{n(n-1)}\sum_{i=1}^n d_{i\pi(i)}
=\frac{n-2}{n(n-1)}U. \enas By symmetry the same is true for the
term $d_{J \pi(J)}$.

Collecting terms and using ${\mathcal F}(U) \subset {\mathcal
F}(\pi)$, where ${\mathcal F}(X)$ denotes the sigma field
generated by the random variable $X$, we have \beas E(U^*-U|U) &=&
\frac{-2}{n(n-1)}(2+2(n-2))U = -\frac{4}{n}U. \enas Thus
(\ref{linearity}) holds with $\lambda=4/n$.

Now we consider the first term in the bound in Theorem \ref{base};
since ${\mathcal F}(U) \subset {\mathcal F}(\pi)$,
\begin{equation} \label{seepi} \mbox{Var}\{E[(U^*-U)^2|U]\} \le
\mbox{Var}\{E[(U^*-U)^2|\pi]\}. \end{equation} From (\ref{wsmw}),
\bea \label{cond-exp} E((U^* - U)^2|\pi) &=& 4E\left( [(d_{IJ} +
d_{\pi(I)\pi(J)})-(d_{I\pi(I)}+d_{J\pi(J)})]^2| \pi \right).\ena

When we expand the square we get the following types of terms; (i)
the square terms from the first group of parentheses, (ii) mixed
terms formed by taking one term from the first group with one term
from the second, (iii) the square terms from the second group,
(iv) mixed terms between values in the first group, and (v) mixed
terms between values in the second group.

(i) The value of the conditional expectation for the square term
E($d\,^2_{IJ} | \pi)$ clearly does not depend on $\pi$, and
therefore contributes a constant value which does not affect the
variance. The same is true for $E(d\,^2_{\pi(I)\pi(J)} | \pi)$
because as $(I,J)$ range over all possible distinct pairs with
equal probability so do $(\pi(I),\pi(J))$.

(ii) Terms such as $E(d_{IJ}d_{I\pi(I)} | \pi)$, evaluate to zero.
In this particular case take expectation over $J$ first and use
$d_{i+}=0$.

By tallying the contributions from terms (iii),(iv), and (v), we
conclude that, up to an additive constant not depending on $\pi$,
and therefore not affecting the variance, (\ref{cond-exp}) equals
\bea \label{withinC}  \frac{8}{n}\sum_{i=1}^n
d_{i\pi(i)}^2 + \frac{8}{n(n-1)}\left( \sum_{|\{i, j\}|=2}
d_{ij}d_{\pi(i)\pi(j)} + \sum_{|\{i,
j\}|=2}d_{i\pi(i)}d_{j\pi(j)}\right). \ena

We may write (\ref{withinC}) as $8\left(A_1 + A_2 +A_3 \right)$
where \beas A_1 = \frac{1}{n}\sum_{i=1}^n d_{i\pi(i)}^2, &\quad&
A_2=
\frac{1}{n(n-1)}\sum_{|\{i, j\}|=2} d_{ij}d_{\pi(i)\pi(j)} \\
\mbox{and} \quad A_3 &=& \frac{1}{n(n-1)}\sum_{|\{i, j\}|=2}
d_{i\pi(i)}d_{j\pi(j)}. \enas

In view of (\ref{seepi}), we now need to compute the variance of
(\ref{withinC}) with respect to a uniform $\pi \in \Pi_n$. We have
$$ \mbox{Var}(8(A_1+A_2+A_3)) \le 8^2 \cdot 3\left(
\mbox{Var}(A_1)+\mbox{Var}(A_2)+\mbox{Var}(A_3)\right).$$

To calculate $\mbox{Var}(A_1)$, apply Lemma \ref{effigy} with
$g_{ij}=d_{ij}^2$ to obtain
$$
\mbox{Var}(A_1)= \frac{1}{n^2(n-1)(n-3)} \left((2n-5)
\sum_{|\{i,j\}|=2}f_{ij}^2 + \sum_{|\{i,j\}|=2}f_{ij}f_{ji}
\right).
$$

For the second term above, by Cauchy-Schwarz, \bea
\label{asyosisays} |\sum_{|\{i, j\}|=2} f_{ij}f_{ji} | \le \sqrt{
\sum_{|\{i, j\}|=2}f^2_{ij} \sum_{|\{i, j\}|=2}f^2_{ji}} =
\sum_{|\{i, j\}|=2}f^2_{ij}. \ena

Collecting terms we conclude \beas \mbox{Var}(A_1) &\le&
\frac{2(n-2)}{n^2(n-1)(n-3)} \sum_{|\{i,j\}|=2}f_{ij}^2 \le
\frac{3}{n^3} \sum_{|\{i,j\}|=2}d_{ij}^4 \enas for $n \ge 8$.

We now turn to $\mbox{Var}(A_2)$. With
$$
{\mathcal I}= \{ {\bf I}=(i,j,k,l,\pi(i),\pi(j),\pi(k),\pi(l)): i
\not = j, k \not = l, \pi \in \Pi_n \},
$$
it can be shown that when $\pi$ is uniform over $\Pi_n$, the
probability of a given ${\bf I} \in \mathcal I$ satisfying ${|\bf
I}|=s$ is
$$
P({\bf I}) = \frac{1}{[n]_s}, \quad \mbox{$s \in \{2,4,6,8\}$,}
\quad \mbox{where} \quad [n]_s = (n-1)(n-3)\cdots(n-s+1).
$$

For ${\bf I}=(i,j,k,l,i',j',k',l') \in {\mathcal I}$ set $d_{\bf
I}=d_{ij}d_{kl}d_{i'j'}d_{k'l'}$. We then have \bea \label{varA2}
\nonumber \mbox{Var}(\sum_{|\{i, j\}|=2} d_{ij}d_{\pi(i)\pi(j)})
&\le& \sum_{|\{i, j\}|=2 |\{k, l\}|=2} d_{ij}d_{kl}
E(d_{\pi(i)\pi(j)}d_{\pi(k)\pi(l)}) \\ &=& \sum_{{\bf I} \in
{\mathcal I}}\frac{1}{[n]_{|{\bf I}|}} d_{\bf I} = \sum_{s \in
\{2,4,6,8\}} \frac{1}{[n]_{s}} \sum_{{\bf I} \in {\mathcal I}(s)}
d_{\bf I}, \ena where ${\mathcal I}(s)$ are all those ${\bf I} \in
{\mathcal I}$ with $|{\bf I}|=s$.

Consider first the case of $s=8$. Since $d_{k'+}=0$, summing over
$l' \not \in \{i,j,\,k,\,l,\,i',j'\}$ we have \bea \label{ok8} &&
\sum_{{\bf I} \in {\mathcal I}(8)} d_{\bf I}= \sum_{{\bf I} \in
{\mathcal I}(8)} d_{ij}d_{kl} d_{i' j'}d_{k'l'} \nonumber \\ &=& -
\sum_{|\{i,j,\,k,\,l,\,i',j',\,k'\}|=7} \sum_{l' \in
\{i,j,\,k,\,l,\,i',j'\}} d_{ij}d_{kl}d_{i'j'}d_{k'l'}. \ena
Applying Cauchy Schwarz to each of the six terms in the inner sum,
the absolute value of the expression is bounded by
$$
6 (n-2)_5 \sum_{|\{i,j\}|=2}d_{ij}^4,
$$
where
$$
(n)_s=n(n-1)\cdots(n-s+1).
$$

For $s \in \{2,4,6\}$ apply Cauchy Schwarz to
$$
\sum_{{\bf I} \in {\mathcal I}(s)} d_{ij} d_{kl} d_{i' j'}
d_{k'l'}
$$
to obtain the bound
$$
(n-2)_{s-2} \sum_{|\{i,j\}|=2}d_{ij}^4.
$$
Therefore \beas \mbox{Var}(A_2) &\le& \frac{1}{(n(n-1))^2} \left(
\frac{6 (n-2)_5}{[n]_8} + \sum_{s \in \{2,4,6\}}
\frac{(n-2)_{s-2}}{[n]_s} \right)
\sum_{|\{i,j\}|=2}d_{ij}^4 \\
&\le& \frac{7}{n^3} \sum_{|\{i,j\}|=2}d_{ij}^4, \enas where the
latter bound holds for $n \ge 10$ and follows by elementary
calculations.

Although $A_3$ and $A_2$ are not identically distributed, it is
easy to see that the variance of $A_3$ can be bounded in exactly
the same manner.

We obtain from (\ref{withinC}) and the above discussion that \bea
\label{varord} \mbox{Var}\{E[(U^*-U)^2|U]\} \le  (8^2 \cdot 3)\,
\frac{17}{n^3} \sum_{i \ne j} d_{ij}^4. \ena

We now apply Theorem \ref{base} to $W=U/\sigma, \, W^*=
U^*/\sigma$. From (\ref{varU}) we conclude that
$$
\mbox{Var}(U)= \sigma^2 \ge \frac{2}{n }\sum_{|\{i,
j\}|=2}d\,^2_{ij}.
$$
It follows from (\ref{varord}), \bea \label{CC2} &&
\mbox{Var}\{E[(W^*-W)^2|W]\} \le \frac{8^2\cdot 3 \cdot 17}{4
n}\sum_{|\{i, j\}|=2}d_{ij}\,^4/(\sum_{|\{i, j\}|=2}d\,^2_{ij})^2.
\ena

With $$ \alpha=\max|d_{ij}-d_{kl}|, $$ we have $|U^*-U| \le 4
\alpha$, and hence \beas |W^* - W| &\le& \frac{1}{\sigma}|U^*-U|
\le \frac{4 \alpha}{\sigma}\\ & \le & \frac{4 \alpha
\sqrt{n}}{\sqrt{2 \sum_{|\{i, j\}|=2}d\,^2_{ij}}} = A. \enas

Applying Theorem \ref{base} with this $A$, $\lambda = 4/n$ and
using expression (\ref{CC2}), we have \bea \label{th1b} \delta
&\le& 86 n^{1/2} \sqrt{\{ \sum_{|\{i,
j\}|=2}d_{ij}\,^4/(\sum_{|\{i, j\}|=2}d\,^2_{ij})^2 \} } +
\frac{243\alpha^3\,n^{5/2}}{(\sum_{|\{i,
j\}|=2}d\,^2_{ij})^{3/2}}. \nonumber \ena

\section{Acknowledgements}
The authors would like to thank Antonio Forcina for his
penetrating comments.

%\nocite{*}
%\bibliographystyle{metron}
%\bibliography{permutationtest}

\vspace*{0.5in}

\noindent
\begin{minipage}[l]{2.2in}
{\footnotesize
Larry Goldstein\\
Department of Mathematics\\
University of Southern California, DRB 155\\
Los Angeles, CA (USA)\\
larry@math.usc.edu}\\
\end{minipage}
\hfill
\begin{minipage}[l]{2.2in}
{\footnotesize
Yosef Rinott\\
Department of Statistics\\
Hebrew University\\
Jerusalem, Israel\\
rinott@mscc.huji.ac.il}\\
\\
\end{minipage}

\end{document}